\documentclass[12pt,leqno]{article}
\usepackage{amsmath, amsthm, amsfonts, amssymb, color}
\usepackage{mathrsfs}
\theoremstyle{definition}

\newcommand{\scr}[1]{\mathscr #1}
\definecolor{wco}{rgb}{0.5,0.2,0.3}

\numberwithin{equation}{section}

\newcommand{\ua}{\uparrow}

\title{{\bf A remark on global solutions to random 3D vorticity equations for small initial data}\footnote{Supported in
 part by  NSFC (11671035, 11771037) and DFG through  CRC 1283}}
\author{
{\bf    Michael R\"{o}ckner$^{c)}$, Rongchan Zhu$^{a,c)}$, Xiangchan Zhu$^{b,c})$}\thanks{Corresponding author}
\thanks{E-mail address: roeckner@math.uni-bielefeld.de(M.R\"{o}ckner),
zhurongchan@126.com(R.C.Zhu), zhuxiangchan@126.com(X.C.Zhu)}
\\
 \footnotesize{ $^{a)}$Department of Mathematics, Beijing Institute of Technology, Beijing 100081,  China}\\
\footnotesize{ $^{b)}$School of Science, Beijing Jiaotong University, Beijing 100044, China}\\
\footnotesize{  $^{c)}$ Department of Mathematics, University of Bielefeld, D-33615 Bielefeld, Germany}}

\date{}

\begin{document}
\maketitle
\def\R{\mathbb R} \def\EE{\mathbb E} \def\P{\mathbb P}\def\Z{\mathbb Z} \def\ff{\frac} \def\ss{\sqrt}
\def\H{\mathbb H}
\def\HH{\mathbf{H}}
\def\DD{\Delta} \def\vv{\varepsilon} \def\rr{\rho}
\def\<{\langle} \def\>{\rangle} \def\GG{\Gamma} \def\gg{\gamma}
\def\ll{\lambda} \def\LL{\Lambda} \def\nn{\nabla} \def\pp{\partial}
\def\dd{\text{\rm{d}}}
\def\Id{\text{\rm{Id}}}\def\loc{\text{\rm{loc}}} \def\bb{\beta} \def\aa{\alpha} \def\D{\scr D}
\def\E{\scr E} \def\si{\sigma} \def\ess{\text{\rm{ess}}}
\def\beg{\begin} \def\beq{\beg}  \def\F{\scr F}
\def\Ric{\text{\rm{Ric}}}
\def\Vol{\text{\rm{Vol}}}
\def\Var{\text{\rm{Var}}}
\def\Ent{\text{\rm{Ent}}}
\def\Hess{\text{\rm{Hess}}}\def\B{\scr B}
\def\e{\text{\rm{e}}} \def\ua{\underline a} \def\OO{\Omega} \def\b{\mathbf b}
\def\oo{\omega}     \def\tt{\tilde} \def\Ric{\text{\rm{Ric}}}
\def\cut{\text{\rm{cut}}} \def\P{\mathbb P} \def\K{\mathbb K}
\def\ifn{I_n(f^{\bigotimes n})}
\def\fff{f(x_1)\dots f(x_n)} \def\ifm{I_m(g^{\bigotimes m})} \def\ee{\varepsilon}
\def\C{\scr C}
\def\PP{\scr P}
\def\M{\scr M}\def\ll{\lambda}
\def\X{\scr X}
\def\T{\scr T}
\def\A{\mathbf A}
\def\LL{\scr L}\def\LLL{\Lambda}
\def\gap{\mathbf{gap}}
\def\div{\text{\rm div}}
\def\Lip{\text{\rm Lip}}
\def\dist{\text{\rm dist}}
\def\cut{\text{\rm cut}}
\def\supp{\text{\rm supp}}
\def\Cov{\text{\rm Cov}}
\def\Dom{\text{\rm Dom}}
\def\Cap{\text{\rm Cap}}\def\II{{\mathbb I}}\def\beq{\beg{equation}}
\def\sect{\text{\rm sect}}\def\H{\mathbb H}

\begin{abstract}In this paper, we prove that the solution constructed in \cite{BR16} satisfies the stochastic vorticity equations with the stochastic integration being understood in the sense of the integration of controlled rough path introduced in \cite{G04}. As a result, we obtain  the existence and uniqueness of the global solutions to the stochastic vorticity equations in 3D case for the small initial data independent of time, which can be viewed as a stochastic version of the Kato-Fujita result (see \cite{KF62}).
\end{abstract}

\vskip 2cm

\section{Introduction}

Consider the stochastic 3D Navier-Stokes equation on $(0,\infty)\times \mathbb{R}^3$:
\begin{equation}\label{eq1.1} \aligned &dX-\Delta Xdt+(X\cdot \nabla)Xdt=\sum_{i=1}^N(B_i(X)+\lambda_iX)d\beta^i(t)+\nabla \pi dt,
\\&\nabla\cdot X=0,
\\&X(0)=x,\endaligned\end{equation}
where $\{\beta^i\}_{i=1}^N$ is a system of independent Brownian motions on a probability space $(\Omega,\mathcal{F},P)$ with normal filtration $(\mathcal{F}_t)_{t\geq0}$, and $\lambda_i\in\mathbb{R}, x:\Omega\rightarrow \mathbb{R}^3$ is a random variable. Here $\pi$ denotes the pressure, $\Delta$ is the Laplacian on $L^2(\mathbb{R}^3;\mathbb{R}^3)$ and $B_i$ are convolution operators given by
$$B_i(X)(\xi)=\int_{\mathbb{R}^3}h_i(\xi-\bar{\xi})X(\bar{\xi})d\bar{\xi}=(h_i*X)(\xi),\quad \xi\in\mathbb{R}^3,$$
where $h_i\in L^1(\mathbb{R}^3),i=1,...,N$.

Consider the vorticity field
$$ U=\nabla\times X=\textrm{curl} X$$
and apply the curl operator to equation \eqref{eq1.1}. We obtain the transport vorticity equation on $(0,\infty)\times \mathbb{R}^3$:
\begin{equation}\label{eq1.2} \aligned &dU -\Delta Udt+((X\cdot \nabla)U-(U\cdot\nabla)X)dt=\sum_{i=1}^N(h_i*U+\lambda_iU)d\beta^i(t),
\\&U_0(\xi)=(\textrm{curl} x)(\xi),\quad \xi\in\mathbb{R}^3.
\endaligned\end{equation}
The vorticity $U$ is related to the velocity $X$ by the Biot-Savart integral operator
\begin{equation}X_t(\xi)=K(U_t)(\xi)=-\frac{1}{4\pi}\int_{\mathbb{R}^3}\frac{\xi-\bar{\xi}}{|\xi-\bar{\xi}|^3}\times U_t(\bar{\xi})d\bar{\xi}, \quad t\in(0,\infty),\xi\in\mathbb{R}^3.\end{equation}
Then one can rewrite the vorticity equation \eqref{eq1.2} as
 \begin{equation}\label{eq1.4} \aligned &dU -\Delta Udt+((K(U)\cdot \nabla)U-(U\cdot\nabla)K(U))dt=\sum_{i=1}^N(h_i*U+\lambda_iU)d\beta^i(t),
\\&U_0(\xi)=(\textrm{curl} x)(\xi),\quad \xi\in\mathbb{R}^3.
\endaligned\end{equation}

In \cite{BR16} using the transformation
$$U_t=\Gamma_ty_t$$
with $$\Gamma_t=\Pi_{i=1}^N\exp\big(\beta^i_t\tilde{B}_i-\frac{t}{2}\tilde{B}_i^2\big),\quad \tilde{B}_i=B_i+\lambda_iI,$$
the authors transformed \eqref{eq1.4} into the following equation
\begin{equation}\label{eq1.5} \aligned &\frac{dy}{dt}-\Gamma^{-1}_t\Delta (\Gamma_ty_t)dt+\Gamma^{-1}_t((K(\Gamma_ty_t)\cdot \nabla)(\Gamma_ty_t)-(\Gamma_ty_t\cdot\nabla)K(\Gamma_ty_t))=0,
\\&y_0=U_0.\endaligned\end{equation}
In \cite{BR16} the authors proved that if the initial value is small enough (compared to a function depending on the paths of Brownian motions $\beta_i$), then there exists a unique solution $y_t$ (in the mild sense) to \eqref{eq1.5}. However, since the initial value is not $\mathcal{F}_0$-measurable, the process $y_t$ is not $(\mathcal{F}_t)_{t\geq0}$-adapted. Therefore, \eqref{eq1.5} cannot be transformed back into \eqref{eq1.4}.

In this paper we  use the result in \cite{BR16} to construct a global solution to \eqref{eq1.4} for small initial data satisfying the following condition \eqref{eq1.10}. Since $y_t$ is not $(\mathcal{F}_t)_{t\geq0}$-adapted, the corresponding $U_t$ is also not $(\mathcal{F}_t)_{t\geq0}$-adapted. Therefore, the stochastic integral should be understood in the sense of a rough path integral or the Skorohod integral. To use the Skorohod integral and find a solution to \eqref{eq1.4} we have to use the shift operator (see \cite{B91}, \cite{N95}), which destroys the following condition \eqref{eq1.10}. Thus in this paper we  understand the stochastic integral of \eqref{eq1.4} in the sense of a rough path integral.
\vskip.10in
\textbf{Framework and main result}
\vskip.10in

First we recall the main result in \cite{BR16}. In the following we denote by $L^p,1\leq p\leq\infty$ the space $L^p(\mathbb{R}^3;\mathbb{R}^3)$ with norm $|\cdot|_p$ and by $C_b([0,\infty);L^p)$ the space of all bounded and continuous functions $u:[0,\infty)\rightarrow L^p$ with the sup norm. We also set $D_i=\frac{\partial}{\partial\xi_i},i=1,2,3$. We set for $p\in (\frac{3}{2},3),q\in(1,\infty)$
$$\eta_t=\|\Gamma_t\|_{L(L^p,L^p)}\|\Gamma_t\|_{L(L^{\frac{3p}{3-p}},L^{\frac{3p}{3-p}})}\|\Gamma^{-1}_t\|_{L(L^q,L^q)}, \quad t\geq0,$$
where $\|\cdot\|_{L(L^p,L^p)}$ is the norm of the space $L(L^p,L^p)$ of linear continuous operators on $L^p.$

For $p\in [1,\infty)$ we denote by $\mathcal{Z}_p$ the space of all functions $y:[0,1]\times\mathbb{R}^3\rightarrow\mathbb{R}^3$ such that
$$t^{1-\frac{3}{2p}}y_t\in C_b([0,\infty);L^p),$$
$$t^{\frac{3}{2}(1-\frac{1}{p})}D_iy_t\in C_b([0,\infty);L^p),\quad i=1,2,3.$$
The space $\mathcal{Z}_p$ is endowed with the norm
$$\|y\|=\sup\{t^{1-\frac{3}{2p}}|y_t|_p+t^{\frac{3}{2}(1-\frac{1}{p})}|D_iy_t|_p;t\in(0,\infty),i=1,2,3\}.$$

In the following we take $\lambda_i\in\mathbb{R}$ such that
$$|\lambda_i|>(\sqrt{12}+3)|h_i|_1, \quad i=1,2,...,N.$$ Consider the equation \eqref{eq1.5} in the following mild sense:
\begin{equation}\label{eq1.6} \aligned &y_t=e^{t\Delta}U_0+\int_0^te^{(t-s)\Delta}\Gamma^{-1}_sM(\Gamma_sy_s)ds, \quad t\in(0,\infty),\endaligned\end{equation}
where $$M(u)=-(K(u)\cdot \nabla)(u)+(u\cdot\nabla)K(u).$$
 The following is the main result in \cite{BR16}.

\beg{thm}\label{T1.1} Let $p,q\in(1,\infty)$ such that
$$\frac{3}{2}<p<2, \frac{1}{q}=\frac{2}{p}-\frac{1}{3}.$$ Let $\Omega_0=\{\sup_{t\geq0}\eta_t<\infty\}$ and consider \eqref{eq1.6} for fixed $\omega\in \Omega_0$. Then $P(\Omega_0)=1$ and there exists a positive constant $C^*$ independent of $\omega\in \Omega_0$ such that, if $U_0\in L^{3/2}$ satisfying
\begin{equation}\label{eq1.10}\sup_{t\geq0}\eta_t|U_0|_{3/2}\leq C^*,\end{equation}
then there exists a unique solution $y\in \mathcal{Z}_p$ to \eqref{eq1.6}. Moreover, for each $\varphi\in L^3\cap L^{\frac{q}{q-1}}$, the function $t\rightarrow \int_{\mathbb{R}^3}y(t,\xi)\varphi (\xi)d\xi$ is continuous on $[0,\infty).$
\end{thm}

To formulate our first main result we introduce the following notations and definitions from rough paths theory: Fix $\frac{1}{3}<\alpha<\frac{1}{2}, 0\leq s<t$, for $X\in C([s,t],\mathbb{R}^N)$  we define
$$\delta X_{uv}:=X_v-X_u,\quad \|X\|_{\alpha,[s,t]}:=\sup_{u,v\in[s,t],u\neq v}\frac{|\delta X_{uv}|}{|u-v|^\alpha}.$$
Moreover, for a tensor process $\mathbb{X}\in C([s,t]^2,\mathbb{R}^{N\times N})$ we define
$$\|\mathbb{X}\|_{2\alpha,[s,t]}:=\sup_{u,v\in[s,t],u\neq v}\frac{|\mathbb{X}_{uv}|}{|u-v|^{2\alpha}}.$$
In fact, $(X,\mathbb{X})$ is an $\alpha$-H\"{o}lder rough path in the sense of \cite{FH14}, Def.2.1 if $\|X\|_{\alpha,[s,t]}<\infty, \|\mathbb{X}\|_{2\alpha,[s,t]}<\infty$ and the following holds for every triple of times $(u,v,w)$
$$\mathbb{X}_{uv}-\mathbb{X}_{uw}-\mathbb{X}_{wv}=\delta X_{uw}\otimes \delta X_{wv}.$$ For an $N$-dimensional Brownian motion $\beta$ on the probability space $(\Omega,\mathcal{F},P)$ and $\mathbb{B}_{uv}:=\int_u^v\delta \beta_{ur}\otimes d\beta_r\in \mathbb{R}^{N\times N}$, it is well known that there exists a set $\Omega_1$ with $P(\Omega_1)=1$ such that for $\omega\in\Omega_1$ $(\beta(\omega),\mathbb{B}(\omega))$ is  an $\alpha$-H\"{o}lder rough path (see \cite{FH14}, Prop. 3.4), where the stochastic integration is understood in the sense of It\^{o}. In the following we consider the problem on $\Omega_1$ $\omega$-wise.
We also introduce the following smaller space for later use: for $\varepsilon>0$
we set $$\mathcal{Z}^\varepsilon_p:=\{y\in\mathcal{Z}_p|\sup_{s\leq u<v\leq t}u^{2\varepsilon+1-\frac{3}{2p}}\frac{|\delta y_{uv}|_p}{|u-v|^\varepsilon}+u^{2\varepsilon+\frac{3}{2}-\frac{3}{2p}}\frac{\sum_{j=1}^3|\delta (D_j y)_{uv}|_p}{|u-v|^\varepsilon}<\infty,\quad 0<s<t\}.$$
Now we recall the notion of a controlled path $Y$ relative to some reference path $X$ due to Gubinelli \cite{G04}.

\beg{defn} Given a  path $X\in C^\alpha([s,t],\mathbb{R}^N)$, we say that $Y\in C^\alpha([s,t],\mathbb{R}^N)$ is controlled by $X$ if there exists $Y'\in C^\alpha([s,t],\mathbb{R}^{N\times N})$ so that the remainder term $R$, for $s\leq u<v\leq t$ given by the formula
$$\delta Y^\mu_{uv}=\sum_{\nu=1}^NY'^{\mu\nu}_u\delta X^\nu_{uv}+R^{\mu}_{uv},$$
satisfies $\|R\|_{2\alpha,[s,t]}<\infty.$  Here the super-index relates to the coordinate.
\end{defn}

By \cite{G04}, if we are given a path  $Y$ controlled by $X$, then we can define the integration of $Y$ against $(X,\mathbb{X})$, which is an extension of Young's integral (see Theorem 1 and Corollary 2 in \cite{G04}): for $0\leq s<t\leq T$
\begin{equation}\label{eqdef}\int_s^tY^\mu dX^\nu:=\lim_{|\mathcal{P}|\rightarrow0}\sum_{i=0}^{n-1}(Y^\mu_{t_i}\delta X^\nu_{t_it_{i+1}}+\sum_{\mu'=1}^NY'^{\mu\mu'}_{t_i}\mathbb{X}^{\mu'\nu}_{t_it_{i+1}}),\end{equation}
where $\mathcal{P}=\{t_0,t_1,...,t_n\}$ is a partition of the interval $[s,t]$ such that $t_0=s,t_n=t,t_{i+1}>t_i,|\mathcal{P}|=\sup_i|t_{i+1}-t_i|$.

Now we give the definition of solutions to equation \eqref{eq1.4}. In the following we define the analytic weak solution to equation \eqref{eq1.4} and we use $\langle \cdot,\cdot\rangle$ to denote the $L^2$ inner product.

\beg{defn}
We say that $U$ is a solution to equation \eqref{eq1.4} if $\Gamma^{-1}U\in\mathcal{Z}^\varepsilon_p$ for some $\varepsilon>0$ and for any $\varphi\in C_c^\infty(\mathbb{R}^3;\mathbb{R}^3),$ the function $t\rightarrow \langle \Gamma^{-1}_tU_t,\varphi\rangle$ is continuous on $[0,\infty)$ and for $0<s<t$,
\begin{equation}\label{eq1.7}\aligned \langle U_t-U_s,\varphi\rangle-\int_s^t[\langle U_r,\Delta\varphi\rangle-\langle M(U_r),\varphi\rangle ]dr=\sum_{i=1}^N\int_s^t\langle \tilde{B}_iU_r,\varphi\rangle d\beta^i_r,\endaligned\end{equation}
$$ U|_{t=0}=U_0,$$
where the integral $\int_s^t\langle \tilde{B}_iU_r,\varphi\rangle d\beta^i_r$ is understood in the sense of \eqref{eqdef} with respect to the rough paths $(\beta, \mathbb{B})$. Here for $0<s<t$ $\langle \tilde{B}_iU,\varphi\rangle\in C^\alpha([s,t])$ is controlled by $\beta$ in the sense of Definition 1.1 and
\begin{equation}\label{eq1.8}\delta(\langle \tilde{B}_iU,\varphi\rangle)_{st}=\sum_{k=1}^N\langle \tilde{B}_k\tilde{B}_i U_s,\varphi\rangle\delta \beta^k_{st}+R^i_{st},\end{equation}
with  $R$ being the remainder term satisfying
 \begin{equation}\label{eq1.9}\|\langle \tilde{B}_k\tilde{B}_i U,\varphi\rangle\|_{\alpha,[s,t]}<\infty,\quad\|R^i\|_{2\alpha,[s,t]}<\infty.\end{equation}
\end{defn}

\beg{Remark} (i) Here due to the singularity of solution $U$ at $t=0$, the stochastic integral defined in \eqref{eqdef} has some problem at $t=0$. So, in \eqref{eq1.7} we only assume $0<s<t$. Since $\Gamma^{-1}U\in\mathcal{Z}_p$, $\int_s^t\langle M(U_r),\varphi\rangle dr$ is well-defined due to (2.35) in \cite{BR16}.

(ii) In general rough paths theory, often approximations are used to give a meaning to the solution of  stochastic equations (see \cite{FH14}, Chapter 12). However, in this case if we need the approximation equations to be well-posed  for small initial data, then  the conditions on the initial value might be artificial. Therefore, since our aim is to prove a stochastic version of the Kato-Fujita result (see \cite{KF62}), the above definition is more suitable. We also want to mention that such kind of definition has also been used for the linear equation in \cite{DFS14}.
\end{Remark}

The  main result of this paper is the following theorem:

\beg{thm}\label{T1.2}  Under the condition of Theorem 1.1 and for $y$ as obtained in Theorem 1.1,
for $\omega\in \Omega_0\cap \Omega_1$, $U_t(\omega):=\Gamma_t(\omega)y_t(\omega)$ is the unique solution to \eqref{eq1.4} in the sense of Definition 1.2.
\end{thm}

\section{Proof of Theorem 1.3}
First, we  prove the following lemma.

\beg{lem}\label{L2.1} ( mild solution $\Leftrightarrow$ weak solution) If $y\in \mathcal{Z}_p$  is the unique solution to \eqref{eq1.6}, then for any $\varphi\in C_c^\infty(\mathbb{R}^3;\mathbb{R}^3)$
\begin{equation}\label{eq2.3} \aligned \langle y_t,\varphi\rangle=&\langle U_0,\varphi\rangle +\int_0^t\big[ \langle y_s,\Delta\varphi\rangle
+\langle \Gamma^{-1}_sM(\Gamma_sy_s),\varphi\rangle\big]ds,\quad t\in[0,\infty).
\endaligned\end{equation}
Conversely, if there exists $y\in \mathcal{Z}_p$ satisfying equation \eqref{eq2.3} for any $\varphi\in C_c^\infty(\mathbb{R}^3;\mathbb{R}^3)$, then $y$ is a solution to \eqref{eq1.6}.
\end{lem}

\proof mild solution $\Rightarrow$ weak solution:  By \eqref{eq1.6} we know that for $\varphi\in C_c^\infty(\mathbb{R}^3;\mathbb{R}^3)$
$$ \aligned \int_0^T\langle y_t,\Delta \varphi\rangle dt=&\int_0^T\langle e^{t\Delta}U_0, \Delta \varphi\rangle dt \\&+\int_0^T\langle\int_0^te^{(t-s)\Delta}\Gamma^{-1}_sM(\Gamma_sy_s)ds,\Delta \varphi\rangle dt.
\endaligned$$
Following similar arguments as in the proof of \cite{PR07}, Proposition G.0.9, we have
$$\int_0^T\langle e^{t\Delta}U_0, \Delta \varphi\rangle dt=\int_0^T\langle U_0, \frac{d}{dt}e^{t\Delta}\varphi\rangle dt=\langle e^{T\Delta}U_0,\varphi \rangle-\langle U_0,\varphi\rangle.$$
$$\aligned &\int_0^T\langle\int_0^te^{(t-s)\Delta}\Gamma^{-1}_sM(\Gamma_sy_s)ds,\Delta \varphi\rangle dt=\int_0^T\langle\Gamma^{-1}_sM(\Gamma_sy_s),(e^{(T-s)\Delta}-I) \varphi\rangle ds.\endaligned$$
Combining the above arguments we have
$$ \aligned \int_0^t\langle y_s,\Delta \varphi\rangle ds=&\langle e^{t\Delta}U_0,\varphi \rangle-\langle U_0,\varphi \rangle+\int_0^t\langle e^{(t-s)\Delta}\Gamma^{-1}_sM(\Gamma_sy_s), \varphi\rangle ds\\&-\int_0^t\langle \Gamma^{-1}_sM(\Gamma_sy_s), \varphi\rangle ds,\endaligned$$
which implies \eqref{eq2.3}.

weak solution $\Rightarrow$ mild solution: By \eqref{eq2.3} and similar arguments as in the proof of \cite{PR07}, Lemma G.0.10, we have for $\zeta\in C^1([0,T];C_c^\infty(\mathbb{R}^3;\mathbb{R}^3))$
\begin{equation}\label{eq2.1}\aligned \langle y_t,\zeta_t\rangle=&\langle U_0,\zeta_0\rangle +\int_0^t\big[ \langle y_s,\Delta\zeta_s+\zeta_s'\rangle
+\langle \Gamma^{-1}_sM(\Gamma_sy_s),\zeta_s\rangle\big]ds,\quad t\in[0,\infty).
\endaligned\end{equation}
Choosing $\zeta_s:=e^{(t-s)\Delta}\varphi$, $\varphi\in C_c^\infty(\mathbb{R}^3;\mathbb{R}^3)$, we have
$$ \aligned \langle y_t,\varphi\rangle=&\langle U_0,e^{t\Delta}\varphi\rangle +\int_0^t\langle e^{(t-s)\Delta}\Gamma^{-1}_sM(\Gamma_sy_s),\varphi\rangle ds.
\endaligned$$
Thus \eqref{eq1.6} follows.  $\hfill\Box$
\vskip.10in
Now we prove the following estimate for the solutions:

\beg{lem}\label{L2.2}  For $T>0,\varphi\in L^{q/(q-1)}\cap L^3$, on $\Omega_0$ $\sup_{t\in[0,T]}|\langle \Gamma_ty_t,\varphi\rangle|<\infty$ and $y\in\mathcal{Z}^\varepsilon_p$ for $0<\varepsilon<\frac{1}{2}-\frac{3}{4p}$, with $p,q$ as in Theorem 1.1.
\end{lem}

\proof We have$$y_t=e^{t\Delta}U_0+\int_0^te^{(t-s)\Delta}\Gamma^{-1}_sM(\Gamma_sy_s)ds.$$
Then on $\Omega_0$
$$\aligned|\langle \Gamma_ty_t,\varphi\rangle|\leq& C\|\Gamma_t\|_{L(L^{3/2},L^{3/2})}|e^{t\Delta}U_0|_{3/2}+C\|\Gamma_t\|_{L(L^{q},L^{q})}\int_0^t|\Gamma^{-1}_sM(\Gamma_sy_s)|_qds
\\\leq& C\|\Gamma_t\|_{L(L^{3/2},L^{3/2})}|U_0|_{3/2}+C\|\Gamma_t\|_{L(L^{q},L^{q})}\int_0^t\|\Gamma^{-1}_s\|_{L(L^q,L^q)}|M(\Gamma_sy_s)|_qds
\\\leq& C\|\Gamma_t\|_{L(L^{3/2},L^{3/2})}|U_0|_{3/2}+C\|\Gamma_t\|_{L(L^{q},L^{q})}\|y\|^2\sup_{s\in[0,t]} \eta_s\int_0^ts^{-5/2+3/p}ds
\\<&\infty,\endaligned$$
where in the second inequality we used (2.15) in \cite{BR16} and in the third inequality we used (2.35) in \cite{BR16} and in the last inequality we used that $\|y\|\leq C|U_0|_{3/2}$ by the proof of Theorem 1 in \cite{BR16}.
Now we prove $y\in\mathcal{Z}_p^\varepsilon$. We have
$$\aligned |\delta y_{uv}|_p\leq &|(e^{v\Delta}-e^{u\Delta})U_0|_p+|(e^{(v-u)\Delta}-1)\int_0^ue^{(u-s)\Delta}\Gamma_s^{-1}M(\Gamma_sy_s)ds|_p
\\&+|\int_u^ve^{(v-s)\Delta}\Gamma_s^{-1}M(\Gamma_sy_s)ds|_p.\endaligned$$
For the first term we have
$$\aligned&|(e^{v\Delta}-e^{u\Delta})U_0|_p=|(e^{(v-u)\Delta}-I)e^{u\Delta}U_0|_p\leq C|(e^{(v-u)\Delta}-I)e^{u\Delta}U_0|_{B^\varepsilon_{p,\infty}}
\\\leq &C(v-u)^\varepsilon |e^{u\Delta}U_0|_{B^{3\varepsilon}_{p,\infty}}\leq C(v-u)^\varepsilon u^{-2\varepsilon}|e^{u\Delta/2}U_0|_{p}\leq C(v-u)^\varepsilon u^{-2\varepsilon-1+\frac{3}{2p}}|U_0|_{3/2},\endaligned$$
where $B^s_{m,n}$ is the usual Besov space and we used the properties of Besov spaces (see \cite{BCD11, GIP15}).
For the second term similarly we have
$$\aligned&|(e^{(v-u)\Delta}-1)\int_0^ue^{(u-s)\Delta}\Gamma_s^{-1}M(\Gamma_sy_s)ds|_p
\\\leq &C(v-u)^\varepsilon \int_0^u|e^{(u-s)\Delta}\Gamma_s^{-1}M(\Gamma_sy_s)|_{B^{3\varepsilon}_{p,\infty}}ds
\\\leq &C(v-u)^\varepsilon \int_0^u(u-s)^{-2\varepsilon}|e^{(u-s)\Delta/2}\Gamma_s^{-1}M(\Gamma_sy_s)|_{p}ds
\\\leq &C(v-u)^\varepsilon \sup\eta_s\|y\|^2\int_0^u(u-s)^{-2\varepsilon-\frac{1}{2}(\frac{3}{p}-1)}s^{-\frac{5}{2}+\frac{3}{p}}ds
\\\leq &C(v-u)^\varepsilon u^{-1-2\varepsilon+\frac{3}{2p}}\sup\eta_s\|y\|^2,\endaligned$$
where in the third inequality we used a similar calculation as (2.17) in \cite{BR16}.
For the third term we have
$$\aligned &|\int_u^ve^{(v-s)\Delta}\Gamma_s^{-1}M(\Gamma_sy_s)ds|_p
\\\leq & C\sup\eta_s\|y\|^2\int_u^v(v-s)^{-\frac{1}{2}(\frac{3}{p}-1)}s^{-\frac{5}{2}+\frac{3}{p}}ds
\\=& C\sup\eta_s\|y\|^2(v-u)^{\frac{3}{2}-\frac{3}{2p}}\int_0^1(1-l)^{-\frac{1}{2}(\frac{3}{p}-1)}[u+l(v-u)]^{-\frac{5}{2}+\frac{3}{p}}dl
\\\leq &C\sup\eta_s\|y\|^2(v-u)^{2\varepsilon}u^{-1-2\varepsilon+\frac{3}{2p}}\int_0^1(1-l)^{-\frac{1}{2}(\frac{3}{p}-1)}l^{-\frac{3}{2}+\frac{3}{2p}+2\varepsilon}dl,\endaligned$$
where we used interpolation in the last inequality. Combining the argument above we obtain that
$$\aligned |\delta y_{uv}|_p\leq &C(v-u)^\varepsilon u^{-2\varepsilon-1+\frac{3}{2p}}(|U_0|_{3/2}+ \sup\eta_s\|y\|^2).\endaligned$$
Similarly we have
$$\aligned |\delta (D_jy)_{uv}|_p\leq &|(e^{v\Delta}-e^{u\Delta})D_jU_0|_p+|(e^{(v-u)\Delta}-1)\int_0^ue^{(u-s)\Delta}D_j\Gamma_s^{-1}M(\Gamma_sy_s)ds|_p
\\&+|\int_u^ve^{(v-s)\Delta}D_j\Gamma_s^{-1}M(\Gamma_sy_s)ds|_p
\\ \leq &C(v-u)^\varepsilon u^{-2\varepsilon-\frac{3}{2}+\frac{3}{2p}}(|U_0|_{3/2}+ \sup\eta_s\|y\|^2),\endaligned$$
where we used a similar calculation as (2.18) in \cite{BR16}. Thus the second result follows.
$\hfill\Box$
\vskip.10in
\emph{Proof of Theorem 1.3}[Existence] Now we check that $U=\Gamma y$ satisfies equation \eqref{eq1.7}. We first calculate $\langle (\delta\Gamma y)_{uv},\varphi\rangle$: for $0<u<v$
$$\aligned\langle (\delta\Gamma y)_{uv},\varphi\rangle=&\langle \delta\Gamma_{uv}y_u,\varphi\rangle+\langle \Gamma_u\delta y_{uv},\varphi\rangle+\langle\delta\Gamma_{uv} \delta y_{uv},\varphi\rangle
\\:=&I_1+I_2+I_3.\endaligned$$
Since $\Gamma_u\varphi=\Pi_{i=1}^N\exp\big(\beta^i_u\tilde{B}_i-\frac{u}{2}\tilde{B}_i^2\big)\varphi$ for $\varphi\in C_c^\infty(\mathbb{R}^3;\mathbb{R}^3)$,
by Taylor expansion  we have
$$\delta\Gamma_{uv}\varphi=\Gamma_u\sum_{i=1}^N(\delta\beta^i_{uv}\tilde{B}_i\varphi-\frac{(v-u)}{2}\tilde{B}_i^2\varphi+\sum_{k=1}^N\frac{1}{2}\tilde{B}_i
\tilde{B}_k\varphi\delta\beta^k_{uv}\delta\beta^i_{uv})+o(|v-u|).$$
Here and in the following $o(|u-v|)$ means a higher order term of $|u-v|$.
Now we recall the following result from Section 3.3 in \cite{FH14}: \begin{equation}\label{eq2.2}\mathbb{B}^{ik}_{uv}+\frac{1}{2}\delta^{ik}(v-u)=\mathbb{B}^{ik}_{str,uv},\end{equation} \begin{equation}\label{eq2.4}\frac{1}{2}(\mathbb{B}^{ik}_{str,uv}+\mathbb{B}^{ki}_{str,uv})=\frac{1}{2}\delta\beta_{uv}^i\delta\beta_{uv}^k,\end{equation}  where $\delta^{ik}=1$ if $i=k$, zero else, and $\mathbb{B}_{str,uv}:=\int_u^v\delta \beta_{ur}\otimes \hat{d}\beta_r\in \mathbb{R}^{N\times N}$ with the integral  in the Stratonovich sense.
Then by symmetry of $\tilde{B}_i
\tilde{B}_k\varphi$ with respect to $i,k$ we have $$\delta\Gamma_{uv}\varphi=\Gamma_u\sum_{i=1}^N(\delta\beta^i_{uv}\tilde{B}_i\varphi-\frac{(v-u)}{2}\tilde{B}_i^2\varphi+\sum_{k=1}^N\tilde{B}_i
\tilde{B}_k\varphi\mathbb{B}^{ik}_{str,uv})+o(|v-u|),$$
which by \eqref{eq2.2} implies that
$$I_1=\sum_{i=1}^N\langle\Gamma_u \tilde{B}_i y_u,\varphi\rangle\delta \beta^i_{uv}+\sum_{i,k=1}^N\langle \Gamma_u \tilde{B}_k\tilde{B}_i y_u, \varphi\rangle\mathbb{B}^{ki}_{uv}+o(|u-v|).$$
Also since $y$ satisfies equation \eqref{eq2.3} and $y\in \mathcal{Z}_p^\varepsilon$, we have
$$\aligned I_2=&\langle y_u,\Delta \Gamma_u^*\varphi\rangle(v-u)+\langle \Gamma_u^{-1}M(\Gamma_uy_u),\Gamma_u^*\varphi\rangle(v-u)+o(|v-u|)
\\=&\langle \Gamma_u y_u,\Delta \varphi\rangle(v-u)+\langle M(\Gamma_uy_u),\varphi\rangle(v-u)+o(|v-u|),\endaligned$$
where $\Gamma_u^*$ means the dual operator of $\Gamma_u$. Here in the first equality we used the following for $u<s$
\begin{equation}\label{eq2.5}\aligned &|\Gamma_s^{-1}M(\Gamma_sy_s)-\Gamma_u^{-1}M(\Gamma_uy_u)|_q
\\\leq &\|\Gamma_s^{-1}-\Gamma_u^{-1}\|_{L(L^q,L^q)}|M(\Gamma_sy_s)|_q+\|\Gamma_u^{-1}\|_{L(L^q,L^q)}|M(\Gamma_sy_s)-M(\Gamma_uy_u)|_q
\\\leq &C_u|s-u|^\varepsilon,\endaligned\end{equation} where in the last inequality we used a similar calculation as Lemma 2.2 in \cite{BR16}.
By the above calculations we know that
$$I_3=\langle \delta y_{uv},\delta\Gamma_{uv}^*\varphi\rangle=o(|v-u|),$$
where $\delta\Gamma_{uv}^*$ means the dual operator of $\delta\Gamma_{uv}$.
The above calculations and Lemma 2.2 and (2.35) in \cite{BR16} imply that $\langle \tilde{B}_i U,\varphi\rangle$ is controlled by $\beta$ in the sense of Definition 1.1 and satisfies \eqref{eq1.8} and \eqref{eq1.9}. By the above calculations we also obtain that for $0<s<t$
$$\aligned &\langle U_t,\varphi\rangle-\langle U_s,\varphi\rangle
\\=&\sum_{[u,v]\in \mathcal{P}}\langle (\delta\Gamma y)_{uv},\varphi\rangle
\\=&\sum_{[u,v]\in \mathcal{P}} \bigg[\sum_{i=1}^N\langle \Gamma_u \tilde{B}_i y_u, \varphi\rangle\delta \beta^i_{uv}+\sum_{i,k=1}^N\langle \Gamma_u \tilde{B}_k\tilde{B}_i y_u, \varphi\rangle\mathbb{B}^{ki}_{uv}
\\&+\langle \Gamma_u y_u,\Delta \varphi\rangle(v-u)+\langle M(\Gamma_uy_u),\varphi\rangle(v-u)+o(|u-v|)\bigg],\endaligned$$
where $\mathcal{P}$ is a partition of the interval $[s,t]$ similar as above. Taking the limit $|\mathcal{P}|\rightarrow0$, by \eqref{eqdef} we obtain that
 $U=\Gamma y$ satisfies the equation \eqref{eq1.7}.

[Uniqueness] Now we prove the uniqueness of the solution. In fact by Theorem 1.1 we already know that the solution to \eqref{eq1.6} is unique, so we only need to prove that $y=\Gamma^{-1}U$ satisfies \eqref{eq2.3}, which is equivalent to \eqref{eq1.6} by Lemma 2.1. We have for $0<u<v$
$$\aligned\langle \delta(\Gamma^{-1}U)_{uv},\varphi\rangle=&\langle \delta\Gamma^{-1}_{uv}U_u,\varphi\rangle+\langle \Gamma^{-1}_u\delta U_{uv},\varphi\rangle+\langle \delta\Gamma^{-1}_{uv}\delta U_{uv},\varphi\rangle
\\:=&J_1+J_2+J_3.\endaligned$$
Since $\Gamma^{-1}U\in \mathcal{Z}^\varepsilon_p$, we obtain  the H\"{o}lder continuity of $U_u$ when $u>0$. Since  $M(U_u)=M(\Gamma_uy_u)$, then \eqref{eq2.5} implies the H\"{o}lder continuity of $M(U_u)$ when $u>0$.
 Then by Corollary 3 in \cite{G04} we have
$$\aligned J_2=&\langle\delta U_{uv},(\Gamma^{-1}_u)^*\varphi\rangle=\langle y_u,\Delta\varphi\rangle(v-u)+\langle \Gamma_u^{-1}M(\Gamma_uy_u),\varphi\rangle(v-u)
\\&+\sum_{k=1}^N\langle \tilde{B}_k y_u,\varphi\rangle \delta\beta^k_{uv}+\sum_{i,k=1}^N\langle \tilde{B}_i\tilde{B}_k y_u,\varphi\rangle\mathbb{B}^{ik}_{uv}+o(|u-v|),\endaligned$$
where $(\Gamma^{-1}_u)^*$ means the dual operator of $\Gamma^{-1}_u$.
Moreover, since
$$\Gamma_u^{-1}\varphi=\Pi_{i=1}^N\exp(-\beta^i_u\tilde{B}_i+\frac{u}{2}\tilde{B}_i^2)\varphi,$$
by Taylor expansion we have $$\delta\Gamma^{-1}_{uv}\varphi=\Gamma_u^{-1}\sum_{i=1}^N(-\delta\beta^i_{uv}\tilde{B}_i\varphi+\frac{(v-u)}{2}\tilde{B}_i^2\varphi+\sum_{k=1}^N\frac{1}{2}\tilde{B}_i
\tilde{B}_k\varphi\delta\beta^k_{uv}\delta\beta^i_{uv})+o(|v-u|).$$
Thus, we have $$J_1=\langle \sum_{i=1}^N(-\delta\beta^i_{uv}\tilde{B}_iy_u+\frac{(v-u)}{2}\tilde{B}_i^2y_u+\sum_{k=1}^N\frac{1}{2}\tilde{B}_i
\tilde{B}_ky_u\delta\beta^k_{uv}\delta\beta^i_{uv}),\varphi\rangle+o(|v-u|),$$
and $$J_3=\langle\delta U_{uv},(\delta\Gamma^{-1}_{uv})^*\varphi\rangle=-\sum_{k,i=1}^N\langle \tilde{B}_i
\tilde{B}_k y_u,\varphi\rangle \delta\beta^k_{uv}\delta\beta^i_{uv}+o(|u-v|),$$
where $(\delta\Gamma^{-1}_{uv})^*$ means the dual operator of $\delta\Gamma^{-1}_{uv}$.
 Using \eqref{eq2.2}, \eqref{eq2.4}  we obtain that
$$\aligned &\sum_{i,k=1}^N\langle \tilde{B}_i\tilde{B}_k y_u,\varphi\rangle\mathbb{B}^{ik}_{uv}
\\=&\sum_{i,k=1}^N\langle \tilde{B}_i\tilde{B}_k y_u,\varphi\rangle\mathbb{B}^{ik}_{str,uv}-\frac{1}{2}\sum_{i=1}^N\langle \tilde{B}^2_i y_u,\varphi\rangle(v-u)
\\=&\sum_{i,k=1}^N\langle \tilde{B}_i\tilde{B}_k y_u,\varphi\rangle[\frac{\mathbb{B}^{ik}_{str,uv}+\mathbb{B}^{ki}_{str,uv}}{2}+\frac{\mathbb{B}^{ik}_{str,uv}-\mathbb{B}^{ki}_{str,uv}}{2}]-\frac{1}{2}\sum_{i=1}^N\langle \tilde{B}^2_i y_u,\varphi\rangle(v-u)
\\=&\sum_{i,k=1}^N\langle \tilde{B}_i\tilde{B}_k y_u,\varphi\rangle\frac{1}{2}\delta\beta_{uv}^i\delta\beta_{uv}^k-\frac{1}{2}\sum_{i=1}^N\langle \tilde{B}^2_i y_u,\varphi\rangle(v-u).\endaligned$$
Thus, we have that for $0<s<t$
$$\aligned &\langle y_t,\varphi\rangle-\langle y_s,\varphi\rangle
\\=&\sum_{[u,v]\in \mathcal{P}}\langle (\delta\Gamma^{-1} U)_{uv},\varphi\rangle
\\=&\sum_{[u,v]\in \mathcal{P}}\bigg[\langle y_u,\Delta\varphi\rangle(v-u)+\langle \Gamma_u^{-1}M(\Gamma_uy_u),\varphi\rangle(v-u)+o(|u-v|)\bigg],\endaligned$$
where $\mathcal{P}$ is a partition of the interval $[s,t]$ as above. Taking the limit $|\mathcal{P}|\rightarrow0$ we obtain that for $0<s<t$
$$\aligned \langle y_t,\varphi\rangle=&\langle y_s,\varphi\rangle +\int_s^t\big[ \langle y_r,\Delta\varphi\rangle
+\langle \Gamma^{-1}_rM(\Gamma_ry_r),\varphi\rangle\big]dr.
\endaligned$$
 Now letting $s\rightarrow0$, by the continuity of $\langle y_s,\varphi\rangle$ and $y\in \mathcal{Z}_p$ we obtain that $y=\Gamma^{-1}U$ satisfies \eqref{eq2.3}. Thus uniqueness follows.
$\hfill\Box$

\end{document}